\documentclass[12pt]{article}
\def\Box{{\setlength{\unitlength}{1.1 ex}
\begin{picture}(1,1)(-0.2,-0.2)
\put(0,0){\framebox(1,1){}}
\end{picture}}}

\font\elevenbb=msbm10 at 10.95pt

\def\N{\hbox{\elevenbb N}}
\def\M{\hbox{\elevenbb M}}
\def\R{\hbox{\elevenbb R}}

\def\T{\hbox{\elevenbb T}}
\def\Gr{Gr\"obner }
\def \bg #1 {\begin{tabular}{{#1}}}
\def \nd {\end{tabular}}
\newcommand \mwhile {{\bf while}\hspace{0.3cm}}
\newcommand \mrepeat {{\bf repeat}\hspace{0.3cm}}
\newcommand \muntil {{\bf until}\hspace{0.3cm}}
\newcommand \mfore {{\bf for\hspace{0.3cm}each}\hspace{0.3cm}}
\newcommand \mdo {{\bf do}\hspace{0.3cm}}

\newcommand \mif {{\bf if}\hspace{0.3cm}}
\newcommand \mthen {{\bf then}\hspace{0.3cm}}
\newcommand \melse {{\bf else}\hspace{0.3cm}}
\newcommand \mchoose {{\bf choose}\hspace{0.3cm}}
\newcommand \mand {{\bf and}\hspace{0.3cm}}

\newcommand \mbegin {{\bf begin}}
\newcommand \mend {{\bf end}}
\newcommand \bb {\hspace{0.3cm}}
\newcommand \h {\hspace{0.5cm}}
\newcommand \hh {\hspace{1.0cm}}
\newcommand \hhh {\hspace{1.5cm}}
\newcommand \hhhh {\hspace{2.0cm}}
\newcommand \hhhhh {\hspace{2.5cm}}
\newcommand \hhhhhh {\hspace{3.0cm}}
\newcommand \hhhhhhh {\hspace{3.5cm}}

\newcounter{cc}
\newcommand \hln {\hfill \addtocounter{cc}{1} \arabic{cc}
            \vskip 0.0cm \noindent }
\newtheorem{definition}{Definition}[section]
\newtheorem{corollary}[definition]{Corollary}
\newtheorem{proposition}[definition]{Proposition}
\newtheorem{example}[definition]{Example}
\newtheorem{theorem}[definition]{Theorem}

\begin{document}
\title{\bf Minimal Involutive Bases
}
\author{Vladimir P. Gerdt \\
       Laboratory of Computing Techniques and Automation\\
       Joint Institute for Nuclear Research\\
       141980 Dubna, Russia \\
       gerdt@jinr.dubna.su
\and
       Yuri A. Blinkov \\
       Department of Mathematics and Mechanics \\
       Saratov University \\
       410071 Saratov, Russia \\
       blinkov@scnit.saratov.su}
\date{}
\maketitle
\begin{abstract}
In this paper we present an algorithm for construction of minimal
involutive polynomial bases which are \Gr bases of the special
form. The most general involutive algorithms are based on the
concept of involutive monomial division which leads to partition of
variables into multiplicative and non-multiplicative. This
partition gives thereby the self-consistent computational procedure
for constructing an involutive basis by performing
non-multiplicative prolongations and multiplicative reductions.
Every specific involutive division generates a particular form of
involutive computational procedure. In addition to three involutive
divisions used by Thomas, Janet and Pommaret for analysis of
partial differential equations we define two new ones. These two
divisions, as well as Thomas division, do not depend on the order
of variables. We prove noetherity, continuity and constructivity of
the new divisions that provides correctness and termination of
involutive algorithms for any finite set of input polynomials and
any admissible monomial ordering. We show that, given an admissible
monomial ordering, a monic minimal involutive basis is uniquely
defined and thereby can be considered as canonical much like the
reduced \Gr basis.
\end{abstract}

\section{Introduction}
\noindent
Computational aspects of constructing \Gr bases
invented by Buchberger~\cite{Buch65} are now under intensive
investigation due to the great theoretical and
practical importance of these bases in computational commutative
algebra and algebraic geometry~\cite{CLO'S95,BWK93,Mishra93}.
\Gr bases are also becoming of greater importance in
non-commutative~\cite{Mora,KRW90,MZ95} and differential
algebra~\cite{Carra87,Ollivier}.

Since its invention about thirty years ago, feasibility of the
Buchberger algorithm has been notably increased. First of all, it was
resulted from discovering criteria for avoiding unnecessary
reductions~\cite{Buch79,Buch85,GM88} which
allow a partial extension to non-commutative case~\cite{MZ95}. Next,
the key role of the reduction and, especially,
selection strategies was experimentally observed, and heuristically
good strategies were found~\cite{Sugar}. For construction of a
lexicographical \Gr basis, which is the most
useful for solving polynomial equations, an efficient computation scheme
was developed in~\cite{FGLM} based on
converting a basis from one ordering into another.

On the other hand, Zharkov and Blinkov~\cite{ZB93} were pioneered in
revealing
another computational scheme for \Gr bases construction in commutative
algebra. They used the partition of variables into multiplicative and
non-multiplicative invented in Pommaret~\cite{Pommaret78} to bring partial
differential equations into so-called involutive form~\cite{Janet} which
has all the integrability conditions satisfied. Zharkov and Blinkov
showed that sequential multiplication of the polynomials in the
system by non-multiplicative variables, and reduction of these prolonged
polynomials modulo others, by means of their multiplicative power products
only,
ends up, under certain conditions, with a \Gr basis. Though the latter
is generally not the reduced basis, it reveals some attractive
features~\cite{Zharkov96}.

Already first computer experiments carried out in~\cite{ZB93} showed
rather high efficiency of the new computational scheme.
However, that algorithm terminates, generally, only for zero-dimensional
ideals and for degree compatible term orderings~\cite{ZB94}.
The algebraic origin of such an algorithmic behavior was analyzed
in~\cite{Apel} where it was also shown that Pommaret involutive bases
are just \Gr ones of ideals in the commutative rings with respect to
non-commutative gradings. Interconnection of Pommaret bases and \Gr bases
was recently investigated also in~\cite{GS95}.

In our previous paper~\cite{GB96} general algorithmic
foundations of involutive approach to commutative algebra were considered,
and a number of new concepts was introduced allowing one to study the
involutive algorithmic procedure in its general form. The central concept
of our analysis is involutive monomial division. Every specific
involutive division generates some particular computation procedure
for constructing the corresponding involutive basis. Every involutive
basis, if it is finite, was proved to be a \Gr basis, generally,
redundant. We formulated the axiomatic properties of an involutive
division which provide a proper partition of variables into multiplicative
and non-multiplicative, and, hence, to construct different divisions.
It was also proved that those partitions used by
Janet~\cite{Janet},
Thomas~\cite{Thomas} and Pommaret~\cite{Pommaret78} are generated by
particular involutive divisions.

Important properties of noetherity, continuity and constructivity
for an involutive division were also characterized. Noetherity provides
for the existence of a finite involutive monomial basis for any monomial
ideal much like to the conventional monomial bases. Continuity
assures involutivity of every locally involutive set. Constructivity
is a strengthening of continuity. It allows one to compute an
involutive monomial basis from the initial one by means of its
enlargement with single non-multiplicative prolongations only, that is,
to avoid
enlargement with multiplicative prolongations. We showed that Janet and
Thomas divisions are noetherian, continuous and constructive whereas
Pommaret division, being continuous and constructive, is not noetherian.
Just by this reason a positive-dimensional polynomial ideal, generally,
does not have a finite Pommaret basis. We presented in~\cite{GB96} a
general form of the involutive algorithm. Its correctness follows from
continuity of a division while termination holds for any
polynomial ideal and for any admissible monomial ordering only
for noetherian divisions. The algorithm involves the Buchberger's
chain criterion to avoid unnecessary reductions.

In the present paper, in addition to Janet, Thomas and Pommaret divisions
analyzed in~\cite{GB96}, we give examples of two more
involutive divisions which are proved to be continuous, constructive
and noetherian. We present also the special form of an involutive
algorithm which, given a constructive noetherian division,
provides computation of a minimal involutive basis. We show that the
monic form of the latter is uniquely defined for any fixed
admissible monomial ordering.

The rest of the paper is organized as follows. In Section 2 we give a
brief review of involutive concepts and methods which are used in the
following sections. In Section 3 we consider some examples of involutive
monomial divisions including those introduced by Thomas, Janet and
Pommaret along with two new ones. In Section 4 we study the minimal
involutive monomial bases. The algorithm for construction of
minimal polynomial bases is described in Section 5, and some concluding
remarks are given in Section 6.

\section{Background of Involutive Approach}

In this section we briefly describe the fundamentals
of the general involutive approach proposed in~\cite{GB96} which
are used in Sections 3-5.

Let $\N$ be a set of non-negative integers, and
$\M=\{x_1^{d_1}\cdots x_n^{d_n}\ |\ d_i\in \N\}$
be a set of monomials in the polynomial ring $\R=K[x_1,\ldots,x_n]$
over zero characteristic field $K$.

By $deg(u)$ and $deg_i(u)$ we denote the total degree of $u\in \M$ and the
degree of variable $x_i$ in $u$, respectively. An admissible monomial
ordering is denoted by $\succ$, and throughout this paper we shall assume
that
\begin{equation}
x_1\succ x_2\succ\cdots\succ x_n\,. \label{var_order}
\end{equation}

The leading monomial and the leading coefficient of polynomial $f\in \R$
with respect to ordering $\prec$ are denoted by $lm(f)$ and $lc(f)$,
respectively. If $F\subset \R$ is a polynomial set, then by $lm(F)$ we
denote the leading monomial set for $F$, and $Id(F)$ will denote the
ideal in $R$ generated by $F$. For the least common multiple and for the
greatest common divisor of two
monomials $u,v\in \M$ we shall use the conventional notations $lcm(u,v)$
and $gcd(u,v)$, respectively.

If monomial $u$ divides monomial $v$ we shall write $u|v$.

\begin{definition}
{\em An {\em involutive division} $L$ on $\M$ is given, if for any finite
monomial set $U\subset \M$ and for any $u\in U$ there is given a
submonoid $L(u,U)$ of $\M$ satisfying the conditions:
\begin{tabbing}
~~(a)~~\=If $w\in L(u,U)$ and $v|w$, then $v\in L(u,U)$. \\
~~(b)  \>If $u,v\in U$ and $uL(u,U)\cap vL(v,U)\not=\emptyset$, then
$u\in vL(v,U)$ or $v\in uL(u,U)$. \\
~~(c)  \> If $v\in U$ and $v\in uL(u,U)$, then $L(v,U)\subseteq L(u,U)$. \\
~~(d)  \> If $V\subseteq U$, then $L(u,U)\subseteq L(u,V)$ for all $u\in V$.
\end{tabbing}
Elements of $L(u,U)$ are called {\em multiplicative} for $u$.
If $w\in uL(u,U)$ we shall write $u|_L w$ and call $u$ {\em ($L-$)involutive
divisor} of $w$. The monomial $w$ in its turn is called {\em ($L-$)involutive
multiple} of $u$. In such an event monomial $v=w/u$ is {\em  multiplicative}
for $u$ and the equality $w=uv$ will be written as $w=u\times v$. If $u$ is
the conventional divisor of $w$ but not involutive one we shall write,
as usual, $w=u\cdot v$. Then $v$ is said to be {\em non-multiplicative} for
$u$.
} \label{inv_div}
\end{definition}

\begin{definition}{\em We shall say that involutive division $L$ is
 {\em globally defined} if for any $u\in \M$ its multiplicative monomials
 are defined irrespective of the monomial set $U\ni u$, that is, if
 $L(u,U)=L(u)$.
} \label{df_global}
\end{definition}

\noindent
Definition~\ref{inv_div} for every $u\in U\subset \M$ provides the partition
\begin{equation}
\{x_1,\ldots,x_n\}=M_L(u,U)\cup NM_L(u,U),\quad M_L(u,U)\cap NM_L(u,U)=\emptyset
\label{part}
\end{equation}
of the set of variables into two subsets: {\em multiplicative}
$M_L(u,U)\subset L(u,U)$ and {\em non-multiplicative}
$NM_L(u,U)\not \in L(u,U)$.
Conversely, if for any finite set $U\in \M$ and any $u\in U$ the
partition~(\ref{part}) is given such that the corresponding
submonoid $L(u,U)$ of monomials in variables in $M_L(u,U)$ satisfies
the conditions (b)-(d), then the partition generates the involutive
division.

The conventional monomial division, obviously, satisfies condition
(b) only in the univariate case.

In what follows monomial sets are assumed to be finite, unless
involutive division $L$ is globally defined. In this case, since
$L$ is defined irrespective to the monomial set, it admits extension
to infinite sets.

\begin{definition}{\em A monomial set $U\in \M$ is {\em involutively
 autoreduced} or {\em $L-$autoreduced} if the condition
 $uL(u,U)\cap vL(v,U)=\emptyset $ holds for all distinct $u,v\in U$.
}\end{definition}

\begin{definition}{\em Given an involutive division $L$, a
 monomial set $U$
 is {\em involutive}\footnote{Janet~\cite{Janet} and
 Thomas~\cite{Thomas} call such sets {\em complete}.}
 with respect to $L$ or $L-$involutive if}
\begin{equation}
 \cup_{u\in U}\,u\,\M = \cup_{u\in U}\,u\,L(u,U)\,.
 \label{imset}
\end{equation}
 \label{inv_mset}
\end{definition}

\begin{definition}{\em An $L-$involutive monomial set $\tilde{U}$
 is called $L-${\em completion} of a set $U\subseteq \tilde{U}$ if
 $$\cup_{u\in U}\,u\,\M = \cup_{u\in \tilde{U}}\,u\,L(u,U)\,.$$
 If there exists a finite $L-$completion $\tilde{U}$ of a finite set
 $U$, then the latter is {\em finitely generated} with respect
 to $L$. The involutive division $L$ is {\em noetherian} if
 every finite set $U$ is finitely generated.
} \label{id_noetherian}
\end{definition}

\begin{proposition}\cite{GB96}
If involutive division $L$ is noetherian, then every monomial ideal
has a finite involutive basis $\bar{U}$.
\label{pr_fin_ib}
\end{proposition}

\begin{proposition}
 If $U$ is a finitely generated monomial set,
 then so is set obtained by autoreduction of $U$ in the sense
 of the conventional monomial division.
 \label{pr_fgs}
\end{proposition}

\noindent
{\bf Proof}\ \ It follows immediately from observation that any
 involutive completion of $U$ is also an involutive completion of its
 autoreduced subset. \hfill{\Box}

\begin{definition}
 {\em A monomial set $U$ is called {\em locally involutive} with respect
 to the involutive division $L$ if
$$
(\forall u\in U)\ (\forall x_i\in NM_L(u,U))\ (\exists v\in U)\
\ [\ v|_L(u\cdot x_i)\ ]\,.
$$
}
\label{mon_loc_inv}
\end{definition}

\begin{definition}
{\em A division $L$ is called {\em continuous} if
for any finite set $U\in \M$ and for any finite
sequence $\{u_i\}_{(1\leq i\leq k)}$ of elements
in $U$ such that
\begin{equation}
(\forall \,i< k)\ (\exists x_j\in NM_L(u_i,U))\ \ [\ u_{i+1}|_L u_i\cdot
x_j\ ] \label{cont_cond}
\end{equation}
the inequality $u_i\neq u_j$ for $i\neq j$ holds.
} \label{def_cont}
\end{definition}

\begin{theorem}\cite{GB96}
If involutive division $L$ is continuous then
local involutivity of any monomial set $U$ implies
its involutivity.
\label{th_cont}
\end{theorem}

\begin{definition}
{\em  A continuous involutive division $L$
 is {\em constructive} if for any $U\subset \M$, $u\in U$,
 $x_i\in NM_L(u,U)$ such that $u\cdot x_i$ has no involutive divisors
 in $U$ and
$$
 (\forall v\in U)\ (\forall x_j\in NM_L(v,U))\ (v\cdot x_j | u\cdot x_i,\
 v\cdot x_j\neq u\cdot x_i)\ \ [\ v\cdot x_j\in
 \cup_{u\in U}\,u\,L(u,U)\ ]
$$
the following condition holds:
\begin{equation}
 (\forall w\in \cup_{u\in U}\,u\,L(u,U))
 \ \ [\ u\cdot x_i\not \in wL(w,U\cup \{w\})\ ].
\label{constr}
\end{equation}
} \label{def_constr}
\end{definition}

\noindent
Given a finite set of polynomials $F\subset \R$
and an admissible ordering $\succ$, multiplicative and
non-multiplicative variables for $f\in F$ are defined in terms of
$lm(f)$ and the leading monomial set $lm(F)$.

The concepts of involutive polynomial reduction and involutive normal form
are introduced similar to their conventional analogues~\cite{Buch85}
with the use of involutive division instead of the conventional one.

\begin{definition} {\em Let $L$ be an involutive division $L$ on $\M$, and
 let $F$ be a finite set of polynomials. Then we shall say:
 \begin{enumerate}
 \renewcommand{\theenumi}{(\roman{enumi})}
 \item $p$ is {\em $L-$reducible} {\em modulo} $f\in F$ if
  $p$ has a term $t=a\,u\in \T$ ($a\neq 0$) such that $u=lm(f)\times v$,
  $v\in L(lm(f),lm(F))$. It yields the {\em $L-$reduction} $p\rightarrow
  g=p-(a/lc(f))\,f\, v$.
 \item $p$ is {\em $L-$reducible modulo} $F$ if there exists $f\in F$ such
  that $p$ is $L-$reducible modulo $f$.
 \item $p$ is {\em in $L-$normal form modulo $F$} if
  $p$ is not $L-$reducible modulo $F$.
\end{enumerate}
} \label{inv_red}
\end{definition}

\noindent
We denote the $L-$normal form of $p$ modulo $F$ by $NF_L(p,F)$. In
contrast, the conventional normal form will be denoted by
$NF(p,F)$. If monomial $u$ is multiplicative for $lm(f)$ ($f\in F$)
and $h=fu$ we shall write $h=f\times u$.

\begin{definition}{\em A finite polynomial set $F$ is {\em $L-$autoreduced}
 if the leading monomial set $lm(F)$ of $F$ is $L-$autoreduced
 and every $f\in F$ does not contain monomials involutively multiple of any
 element in $lm(F)$.
}
\end{definition}

\begin{theorem}\cite{GB96} If set $F\subset \R$ is $L-$autoreduced, then
 $NF_L(p,F)=0$ iff $p\in \R$ is presented in the form
$p=\sum_{ij} c_if_i\times u_{ij}$ where $f_i\in F$, $c_i\in K$, and
$u_{ij}\in L(lm(f),lm(F))$ are such that $u_{ij}\neq u_{ik}$ for $i\neq k$.
\label{th_nfzero}
\end{theorem}

\begin{corollary}\cite{GB96}
If polynomial set $F$ is $L-$autoreduced, then $NF_L(p,F)$ is uniquely
defined
for any $p\in \R$, and $NF_L(p_1+p_2,F)=NF_L(p_1,F)+NF_L(p_2,F)$.
\label{c_inf}
\end{corollary}

\begin{definition}{\em An $L-$autoreduced set $F$ is called
{\em ($L-$)involutive} if
$$
(\forall f\in F)\ (\forall u\in \M)\ \ [\ NF_L(fu,F)=0\ ]\,.
$$
Given $v\in \M$ and an $L-$autoreduced set $F$, if there exist
$f\in F$ such that $lm(f)\prec v$ and
\begin{equation}
(\forall f\in F)\ (\forall u\in \M)\ (lm(f)\cdot u\prec v)\ \
[\ NF_L(fu,F)=0\ ]\,, \label{cond_pinv}
\end{equation}
then $F$ is called {\em partially involutive up to the monomial $v$}
with respect to the admissible ordering $\prec$. $F$ is
still said to be partially involutive up to $v$ if $v\prec lm(f)$ for
all $f\in F$.
} \label{def_inv}
\end{definition}

\begin{theorem}\cite{GB96}
 An $L-$autoreduced set $F\subset \R$ is involutive
 with respect to a continuous involutive
 division $L$ iff the following (local) involutivity conditions
 hold
$$
 (\forall f\in F)\ (\forall x_i\in NM_L(lm(f),lm(F)))\
 \ [\ NF_L(f\cdot x_i,F)=0\ ]\,.
$$
Correspondingly, partial involutivity~(\ref{cond_pinv}) holds
iff
$$
(\forall f\in F)\ (\forall x_i\in NM_L(lm(f),lm(F)))\
(lm(f)\cdot x_i\prec v)\ \ [\ NF_L(f\cdot x_i,F)=0\ ]\,.
$$
\label{th_inv_cond}
\end{theorem}

\begin{theorem}\cite{GB96} If $F\subset \R$ is an $L-$involutive basis, then
 it is also a \Gr basis, and the equality of
the conventional and $L-$normal forms
$NF(p,F)=NF_L(p,F)$ holds for any polynomial $p\in \R$.
\noindent
If set $F$ is partially involutive up to
the monomial $v$, then the equality of the normal forms
$NF(p,F)=NF_L(p,F)$ holds for any $p\in \R$ such that $lm(p)\prec v$.
\label{th_nf}
\end{theorem}

\begin{theorem}\cite{GB96}
Let $F$ be a finite $L-$autoreduced polynomial set, and let
$g\cdot x$ be a non-multiplicative prolongation of $g\in F$.
Then $NF_L(g\cdot x,F)=0$ if the following holds
$$
(\forall h\in F)\ (\forall u\in \M)\ \bigl(\ lm(h)\cdot u\prec
lm(g\cdot x)\  \bigr)\ \
[\ NF_L(h\cdot u,F)=0\ ]\,,
$$
$$
     (\exists f,f_0,g_0\in F)
\left[
\begin{array}{l}
 lm(f_0)|lm(f)\,,\ lm(g_0)|lm(g) \\[0.1cm]
 lm(f)|_L lm(g\cdot x)\,,\ lcm(f_0,g_0)\prec lm(g\cdot x) \\[0.1cm]
 NF_L\bigl(f_0\cdot \frac{lt(f)}{lt(f_0)},F\bigl)=
 NF_L\bigl(g_0\cdot \frac{lt(g)}{lt(g_0)},F\bigl)=0
\end{array}
\right]\,.
\label{inv_crit}
$$
\label{th_criteria}
\end{theorem}

\section{Examples of Involutive Divisions}

First of all, we give three examples of involutive division used
in~\cite{Pommaret78,Janet,Thomas} for
analysis of algebraic differential equations. For the proof of validity of
properties (b)-(d) in Definition~\ref{inv_div} for these divisions we
refer to~\cite{GB96}.

\begin{example}{\em Thomas division~\cite{Thomas}.
Given a finite set $U\subset \M$, the variable $x_i$ is considered as
multiplicative
for $u\in U$ if $deg_i(u)=max\{deg_i(v)\ |\   v\in U\}$, and
non-multiplicative, otherwise.
} \label{div_T}
\end{example}

\begin{example}{\em Janet division~\cite{Janet}. Let set $U\subset \M$
be finite. For each $1\leq i\leq n$ divide $U$ into groups
labeled by non-negative integers $d_1,\ldots,d_i$:
$$[d_1,\ldots,d_i]=\{\ u\ \in U\ |\ d_j=deg_j(u),\ 1\leq j\leq i\ \}.$$
A variable $x_i$ is multiplicative for $u\in U$ if
$i=1$ and $deg_1(u)=max\{deg_1(v)\ |\ v\in U\}$,
or if $i>1$, $u\in [d_1,\ldots,d_{i-1}]$ and
$deg_i(u)=max\{deg_i(v)\ |\ v\in [d_1,\ldots,d_{i-1}]\}$.
} \label{div_J}
\end{example}

\begin{example}{\em
 Pommaret division~\cite{Pommaret78}. For a monomial
$u=x_1^{d_1}\cdots x_k^{d_k}$ with $d_k>0$ the variables $x_j,j\geq k$ are
considered as multiplicative and the other variables as non-multiplicative.
For $u=1$ all the variables are multiplicative.
} \label{div_P}
\end{example}

\noindent
Now we present two more examples of divisions which, as does Thomas
division, do not rest on the variable ordering.

\begin{example} {\em Division I. Let $U$ be a finite monomial set. The
variable
 $x_i$  is non-multiplicative for $u\in U$ if there is $v\in U$ such
 that
 $$ x_{i_1}^{d_1}\cdots x_{i_m}^{d_m}u=lcm(u,v),\quad 1\leq m\leq [n/2],
 \quad d_j>0\ \ (1\leq j\leq m)\,,$$
 and $x_i\in \{x_{i_1},\ldots,x_{i_m}\}$.
} \label{div_I}
\end{example}

\begin{example} {\em Division II. For monomial
 $u=x_1^{d_1}\cdots x_k^{d_n}$ the variable $x_i$ is
 multiplicative if $d_i=d_{max}(u)$ where
 $d_{max}(u)=max\{d_1,\ldots,d_n\}$.
} \label{div_II}
\end{example}

\noindent
To distinguish the above divisions, the related subscripts $T,J,P,I,II$ will
be used.

\noindent
We note that
\begin{itemize}
\item Thomas division, Divisions I and II do not depend on the ordering on
 the
 variables $x_i$. Two other divisions, as defined, are based on the
 ordering~(\ref{var_order}).
\item Pommaret division and Division I are globally
 defined in accordance with Definition~\ref{inv_div}, and, hence,
 admit extension to infinite monomial sets.
\end{itemize}

\begin{proposition}
 Divisions I and II are involutive.
\end{proposition}

\vskip 0.2cm
\noindent
{\bf Proof}\ \ {\em Division I.} First of all, we prove that the
 condition (b) in Definition~\ref{inv_div} is fulfilled.
 Let $u\neq v$ be elements in $U$ such that $u|_Iw$ and $v|_Iw$ for some
 $w\in \M$. If $u|_Iv$ or $v|_Iu$, then we are done. Otherwise, $lcm(u,v)/u$
 or $lcm(u,v)/v$ contains non-multiplicative variables for $u$ or $v$,
 respectively. Because $lcm(u,v)|w$, it follows that $w$ cannot be
 involutively multiple of both $u$ and $v$.

 Consider now $u\in U$ such that $u|_Iv$ for some $v\in U$, and $v\neq u$.
 Suppose $v|_I w$ for some $w\in \M$, and assume for a contradiction
 that $w$ is not involutively multiple of $u$. Then there are variables
 $x_{i_1},\ldots,x_{i_m}$ $(1\leq m \leq [n/2])$ containing in $w/v$
 which are non-multiplicative
 for $u$ and there is $t\in U$ such that
 $ux_{i_1}^{k_1}\cdots x_{i_m}^{k_m} = lcm(u,t)$.
 Because $v/u$ does not contain $x_{i_1},\ldots,x_{i_m}$ it follows
 $vx_{i_1}^{k_1}\cdots x_{i_m}^{k_m} = lcm(v,t)$,
 that contradicts our assumption that $w\in vL(v,U)$ and proves the
 fulfillment of condition (c).

 The condition (d) holds too, since an enlargement of the set $U$ may,
 obviously, only produces extra non-multiplicative variables for any
 $u\in U$.

\vskip 0.1cm
{\em Division II.} Let $u$ with $d_u=d_{max}(u)$ be
 an involutive divisor of some monomial $w\in \M$. Then, by definition,
 $deg_i(u)=min(deg_i(w),d_u)$ $(1\leq i\leq n)$. Thus, given monomial
 $w$ and number $d_u$ such that $d_u\leq d_w$ where $d_w=d_{max}(w)$,
 the corresponding involutive divisor $u$ of $w$
 is uniquely defined. If there are two involutive divisors $u,v$ of $w$
 with $d_u < d_v$, then it follows that
\begin{center}
$ \begin{array}{lcl}
  deg_i(u)=deg_i(v)=deg_i(w) & if & deg_i(w)\leq d_u\,,\\
  d_u < deg_i(v)=min(deg_i(w),d_v) & if & deg_i(w) > d_u\,.
\end{array}$
\end{center}
Hence, $u$ is involutive divisor of $v$ and the
condition (b) is fulfilled.

 The condition (c) is an easy consequence of the relations
 $deg_i(u)=min(deg_i(v),d_u)$ and $deg_i(v)=min(deg_i(w),d_v)$.

 The condition (d) holds trivially, because the division as well as Pommaret
 one does not depend on monomial set $U$ at all.
\hfill{\Box}

\noindent
\begin{proposition} For any finite monomial set $U$ and for any monomial
 $u\in U$, the inclusion $M_T(u,U)\subseteq M_I(u,U)$ and,
 respectively, $NM_I\subseteq NM_T(u,U)$ holds.
 \label{id_T_I}
\end{proposition}
{\bf Proof}\ \
 If $x_i\in NM_I(u,U)$, then, obviously,
 $deg_i(u) < h_i=max\{deg_i(u)\ |\   u\in U\}$,
 and, hence, $x_i\in NM_T(u,U)$.
\hfill{\Box}
\vskip 0.1cm
\noindent
\begin{example} {\em $U=\{x^2,xy,z\}$ ($x\succ y\succ z$).
\begin{center}
\vskip 0.2cm
\bg {|c|c|c|c|c|c|c|c|c|c|c|} \hline\hline
monomial & \multicolumn{2}{c|}{Thomas} &\multicolumn{2}{c|}{Janet}
& \multicolumn{2}{c|}{Pommaret} & \multicolumn{2}{c|}{Division I}
& \multicolumn{2}{c|}{Division II} \\ \cline{2-11}
       & $M_T$ & $NM_T$  &   $M_J$   & $NM_J$   & $M_P$   & $NM_P$
       & $M_I$ & $NM_I$  &   $M_{II}$ & $NM_{II}$ \\ \hline
 $x^2$ & $x$ & $y,z$ & $x,y,z$ & $-$   & $x,y,z$  &  $-$  & $x$ & $y,z$
       & $x$ & $y,z$ \\
 $xy$  & $y$ & $x,z$ & $y,z$   & $x$    & $y,z$ &  $x$  & $y$ & $x,z$
       & $x,y$ & $z$ \\
 $z$   & $z$ & $x,y$ &  $y,z$    & $x$  & $z$   & $x,y$ & $y,z$ & $x$
       & $z$ & $x,y$ \\ \hline \hline
\nd
\vskip 0.2cm
\end{center} \label{exm_1}
}
\end{example}

\vskip 0.1cm
\begin{proposition}~Divisions given by Examples~\ref{div_T}-\ref{div_II}
 are continuous and constructive. All these divisions except that of
 Pommaret are also  noetherian.
\label{th_cont_TJP}
\end{proposition}

\noindent
{\bf Proof}\ \ The proof for Thomas, Janet and Pommaret divisions
is given in~\cite{GB96}. Consider Divisions I and II.

{\sc Continuity}. Let $U$ be a finite set, and
$\{u_i\}_{(1\leq i\leq M)}$
be a sequence of elements in $U$ satisfying the
conditions~(\ref{cont_cond}). In accordance with Definition~\ref{def_cont}
we shall show that there are no coinciding elements in the sequence for
each of the two divisions. There are the following two alternatives:
\begin{equation}
 (i)\ \ u_i=u_{i-1}\cdot x_{j};\qquad (ii)\ \ u_i\neq u_{i-1}\cdot x_{j}\,.
 \label{two_cases}
\end{equation}
Extract from the sequence $\{u_i\}$ the subsequence
$\{t_k\equiv u_{i_k}\}_{(1\leq k\leq K\leq M)}$
of those elements which occur in the left-hand side of
relation $(ii)$ in~(\ref{two_cases}).

{\em Division I.}
Show that $t_k|_I lcm(t_{k-1},t_k)$ and $t_{k}\neq
lcm(t_{k-1},t_{k})$. We have
$t_k\times \tilde{w}_k=u_{i_k-1}\cdot x_{j_k}=t_{k-1}\cdot
\tilde{v}_{k-1}$
where
$\neg \tilde{w}_k|\tilde{v}_{k-1}$. Indeed,
suppose $\tilde{w}_k |\tilde{v}_{k-1}$. Apparently, this implies
the relation
$t_k=u_{l}\cdot z_l$ where
$i_{k-1}\leq l< i_k$, and the variable
$x_{j_l}\in NM_I(u_l,U)$, which figures in Definition~\ref{def_cont} of
the sequence $\{u_i\}$, satisfies $x_{j_l}|\tilde{w}_k$ and
$\neg x_{j_l}|z_l$.
It follows that $lcm(t_k,u_{l+1})=t_{k} x_{j_l}$ what, in
accordance with definition of the division in Example~\ref{div_I},
contradicts multiplicativity of $x_{j_l}$ for $t_k$.
Therefore, we obtain the relation
\begin{equation}
\left\{
\begin{array}{l}
t_k\cdot v_k=t_{k+1}\times w_{k+1}\,, \\
gcd(v_{k},w_{k+1})=gcd(v_k,w_k)=1\,,
\end{array}
\right.
\label{rel_uv}
\end{equation}
where $w_{k+1}$ contains more then $[n/2]$
variables with positive exponents, and, hence,
$v_k$ contains only non-multiplicative variables
for $t_k$.

We claim now that any $v_j$ occurring in~(\ref{rel_uv}) with $j>k$ as well as
$v_k$ contain only non-multiplicative variables for
$t_k$. For $j=k+1$ we multiply $t_k v_k$ by $v_{k+1}$
$$
\left\{
\begin{array}{l}
t_k v_k v_{k+1}=(t_{k+1}\cdot v_{k+1}) w_{k+1}=(t_{k+2}\times w_{k+2})w_{k+1}\,, \\
gcd(v_{k},w_{k+1})=gcd(v_{k+1},w_{k+1})=gcd(v_{k+1},w_{k+2})=1\,.
\end{array}
\right.
$$
It yields
\begin{equation}
\left\{
\begin{array}{l}
t_k \hat{v}_k v_{k+1}=(t_{k+2}\times \hat{w}_{k+2})w_{k+1}\,, \\
gcd(\hat{v}_k v_{k+1},\hat{w}_{k+2}w_{k+1})=1\,.
\end{array}
\right. \label{two-factors}
\end{equation}
Because $w_{k+1}$ contains more than $[n/2]$ variables,
the number of variables occurring in the product $\hat{v}_k v_{k+1}$
is less or equal $[n/2]$, and, thus,
variables which are multiplicative for $t_k$ are not contained in
$v_{k+1}$.

If we proceed, sequentially multiplying the upper equality
in~(\ref{two-factors}) by $v_{k+j}$ ($j=2,\ldots $), rewriting the
right-hand side of every product in terms of $t_{k+j+1}$ and
cancelling the common factors, then we obtain the equality
$$
\left\{
\begin{array}{l}
t_k \hat{v}_k \cdots \hat{v}_{k+j-1} v_{k+j}=
 (t_{k+j+1}\times \hat{w}_{k+j+1})\hat{w}_{k+1}\cdots\hat{w}_{k+j-1}
 w_{k+j}\,, \\
gcd(\hat{v}_k \cdots \hat{v}_{k+j-1}v_{k+j},\hat{w}_{k+j+1}\hat{w}_{k+1}\cdots
 \hat{w}_{k+j-1}w_{k+j})=1\,.
\end{array}
\right.
$$
It proves the claim and implies $t_i\neq t_j$ for $i\neq j$.

It remains to prove that elements of the sequence
$\{u_i\}_{(1\leq i\leq M)}$ which occur in the left-hand side
of relation $(i)$ in~(\ref{two_cases})
are also distinct. Assume for a contradiction that there are two elements
$u_j=u_k$ with $j < k$. In between these elements there is, obviously,
an element from the left-hand side of relation $(ii)$ in~(\ref{two_cases}).
Let
$u_{i_m}$ ($j<i_m<k$) be the nearest such element to $u_j$. Considering
the same non-multiplicative prolongations of $u_k$ as those of $u_j$
in the initial sequence, one can construct a sequence such that the
subsequence of the left-hand sides of relation (ii) in~(\ref{two_cases})
has two identical elements $u_{i_k}=u_{i_m}$ with $i_k>i_m$.

\vskip 0.1cm
{\em Division II.} The above defined elements $t_k$ which occur in the
left-hand side of the relation $(ii)$ in~(\ref{two_cases}) are distinct
because $d_{max}(t_{k+1})<d_{max}(t_k)$. The other elements
occurring in relation $(i)$ in~(\ref{two_cases}) are also distinct since
$deg(u_{i_k+j})=deg(u_{i_k+j-1})+1$ $(j=1,\ldots,i_{k+1}-i_k-1)$ and
$$d_{max}(t_k)=d_{max}(u_{i_{k}+1})=\cdots =d_{max}(u_{i_{k+1}-1})\,.$$

\vskip 0.1cm
{\sc Constructivity}. {\em Division I.} Let $u\cdot x_i$, $u\in U$,
$x_i\in NM_I(u,U)$ be a non-multiplicative prolongation
such that
$$u\cdot x_i=u_1v\times w,\quad u_1\in U,\quad v\in I(u_1,U),\quad
w\in I(u_1v,U\cup \{u_1v\}),\quad w\neq 1\,.$$ Show that if
$x_j|w$, then $x_j\in M_I(u_1,U)$. Suppose $x_j\in NM_I(u_1,U)$. It
means that there is $v_1\in U$ satisfying $deg_j(u_1)<deg_j(v_1)$.
Because $\neg x_j|v$, we have $deg_j(u_1v)<deg_j(v_1)$, and, hence,
$x_j\in NM_I(u_1v,U\cup \{u_1v\})$.

{\em Division II.} Since this division is globally defined, its
constructivity ia an immediate consequence of the property (c) in
Definition~\ref{inv_div}.

\vskip 0.1cm
{\sc Noetherity}. {\em Division I.} Its noetherity follows from
Proposition~\ref{id_T_I} and noetherity of Thomas division, since
every Thomas completion of a set $U$, obviously, is also its
completion with respect to Division I.

{\em Division II.} Given a finite set $U\subset \M$ and $u\in U$ with
$d_u=d_{max}(u)$, complete the set by the monomial
$x_1^{d_u}\cdots x_n^{d_u}$ and all its divisors multiple of $u$. If we
do such a completion for every $u\in U$ we obtain, apparently,
an involutive completion of $U$.
\hfill{\Box}

\section{Minimal Involutive Monomial Bases}

Let $U$ be a finitely generated monomial set with respect to
involutive division $L$. In this case a finite involutive completion
$\tilde{U}\supseteq U$ forms the involutive basis of the monomial ideal
generated by $U$. A monomial ideal may not have the unique
involutively autoreduced basis. For instance,
from the definition of Janet division given in
Example~\ref{div_J} it is easy to see that any finite monomial set is
Janet autoreduced.
Therefore, enlargement of a Janet basis by a prolongation of any
its element and Janet completion of the enlarged set leads to another
Janet basis of the same monomial ideal. Similarly, Thomas division
and Division I do not provide uniqueness of involutively autoreduced bases
whereas Pommaret division and Division II do, as the following
proposition shows.

\begin{proposition} Let $L$ be a globally defined involutive division.
 Then any monomial ideal has the unique $L-$involutive basis.
\label{local_id}
\end{proposition}

\noindent
{\bf Proof}\ \ Assume that there are
two distinct $L-$bases $\bar{U}_1$ and $\bar{U}_2$ of
the monomial ideal $Id(U)$ where $U$ is the finite monomial set
generating the ideal and autoreduced in the sense of the conventional
monomial division.
Both $\bar{U}_1$ and $\bar{U}_2$ are apparently involutive completion
of $U$. It follows $\bar{U}_1\setminus \bar{U}_2\neq \emptyset$ and
$\bar{U}_2\setminus \bar{U}_1\neq \emptyset$.
Otherwise one of sets $\bar{U}_1,\bar{U}_2$ would contain another,
and, hence, could not be $L-$autoreduced. Indeed, let $\bar{U}_2
\subset \bar{U}_1$. Then any element of $u\in \bar{U}_1\setminus \bar{U}_2$
is multiple of some element in $U$, and, in accordance with
Definition~\ref{id_noetherian}, $u$ is involutively multiple
of some element $v\in \bar{U}_2$.

We obtain that for any $u\in \bar{U}_1\setminus \bar{U}_2$ there is
$v\in \bar{U}_2\setminus \bar{U}_1$ such that $v|_Lu$ and
for any $v\in \bar{U}_2\setminus \bar{U}_1$ there is $w\in
\bar{U}_1\setminus \bar{U}_2$ such that $w|_Lv$.
Thus, by property (c) in Definition~\ref{inv_div},
given $u\in \bar{U}_1\setminus \bar{U}_2$ there exist $w \in
\bar{U}_1\setminus \bar{U}_2$ such that $w|_L u$.
Since $\bar{U}_1$ is $L-$autoreduced, it is possible only if $u=w$.
But this implies $u=v$. The obtained contradiction proves the proposition.
\hfill{\Box}

\begin{definition}
{\em Let $L$ be an involutive division, and $Id(U)$ be a monomial ideal. Then
 its $L-$involutive basis $\bar{U}$ will be called {\em minimal} if for any
 other involutive basis $\bar{V}$ of the same ideal the inclusion
 $\bar{U}\subseteq \bar{V}$ holds.
} \label{def_mmi}
\end{definition}

\begin{proposition}
 If $U\subset \M$ is a finitely generated set with respect to a
 constructive involutive division, then monomial ideal $Id(U)$ has
 the minimal involutive basis.
 \label{pr_mmb}
\end{proposition}

\noindent
{\bf Proof}\ \ The proof follows immediately from Proposition~\ref{pr_fgs}
 and existence of the minimal involutive completion for a finitely
 generated set~\cite{GB96}. \hfill{\Box}

\vskip 0.2cm
\noindent
If $L$ is constructive, then to compute the minimal involutive basis
for an ideal generated by a given finite monomial set one can use the
following algorithm which is a slightly modified version of
algorithm {\bf InvolutiveCompletion} in paper~\cite{GB96}.

\setcounter{cc}{00}
\vskip 0.3cm
\noindent
\h Algorithm {\bf MinimalInvolutiveMonomialBasis:}
\vskip 0.2cm
\noindent
\h {\bf Input:}  $U$, a finite monomial set
\vskip 0.0cm
\noindent
\h {\bf Output:} $\bar{U}$, a minimal involutive basis of $Id(U)$
\vskip 0.0cm
\noindent
\h \mbegin
\hln
\hh $\bar{U}:=Autoreduce(U)$
\hln
\hh \mchoose any admissible monomial ordering $\prec$
\hln
\hh \mwhile exist $u\in \bar{U}$ and $x\in NM_L(u,\bar{U})$ s.t.
\hln
\hhh $u\cdot x$ has no involutive divisors in $\bar{U}$\bb \mdo
\hln
\hhh \mchoose such $u,x$ with the lowest $u\cdot x$ w.r.t. $\prec$
\hln
\hhh $\bar{U}:=\bar{U}\cup \{u\cdot x\}$
\hln
\hh \mend
\hln
\h \mend
\hln
\vskip 0.2cm
\noindent
The proof of {\em correctness} and {\em termination}, for a
finitely generated set, of this algorithm is the same as that of
algorithm {\bf InvolutiveCompletion}~\cite{GB96} if
Proposition~\ref{pr_fgs} is taken into account. In effect, the
below algorithm constructs the minimal involutive completion of an
autoreduced, in the sense of the conventional monomial division,
initial monomial set. This autoreduction is just done in line 2 of
the algorithm.

\begin{example} {\em (Continuation of Example~\ref{exm_1}). The
 minimal involutive bases of the ideal generated by
 the set $U=\{x^2,xy,z\}$ ($x\succ y\succ z$) are given by
\begin{eqnarray*}
&\bar{U}_T&=\{x^2,xy,z,xz,yz,x^2y,xyz,x^2z,x^2yz\}\,,\\
&\bar{U}_J&=\{x^2,xy,z,xz\}\,,\\
&\bar{U}_P&=\{x^2,xy,z,xz,yz,y^2z,\ldots,y^kz,\ldots \}\,, \\
&\bar{U}_I&=\{x^2,xy,z,xz,x^2y,xyz,x^2z,x^2yz\}\,,\\
&\bar{U}_{II}&=\{x^2,xy,z,xz,yz,xyz\}\,,
\end{eqnarray*}
where $k\in \N$ ($k>2$), and subscripts in the left-hand sides
stand for different involutive divisions considered in Section 3.
This example explicitly shows that Pommaret division is not noetherian.
However, for
another ordering $z\succ x\succ y$ the set $U$ is finitely generated, and
then $\bar{U}_P=U$.} \label{exm_2}
\end{example}

One should note that selection of a $L-$irreducible non-multiplicative
prolongation which is lowest with respect to an admissible
monomial ordering and which we call {\em normal} is of fundamental
importance for the above algorithm. We demonstrate this fact by the
following example.

\begin{example} {\em Let $U=\{x^2,xz,y\}$ and $L$ be Pommaret division with
 $x\succ y\succ z$. By the normal selection strategy, the lowest irreducible
 non-multiplicative prolongation is $y\cdot x$ with respect to any admissible
 monomial ordering. Enlargement of $U$ by $xy$ gives the Pommaret
 basis $\bar{U}=\{x^2,xy,xz,y\}$ of ideal $Id(U)$ which is obviously minimal.
 This shows that $U$ is a finitely generated set.
 However, if we would take first the prolongation $xz\cdot y$ which is
 involutively irreducible modulo $U$, but not lowest, then we might obtain
 the infinite chain of irreducible prolongations:
 $$ xz\rightarrow xyz\rightarrow xy^2z \cdots \rightarrow xy^kz \rightarrow
 \cdots $$
 } \label{exmp_nmp}
\end{example}

\begin{definition}
{\em Let $L$ be a constructive involutive division, $U$ be a finite monomial
 set and $V=Autoreduce(U)$. Then set $U$ will be called {\em $(L-)$compact}
 if $U=V$ or $U$ is obtained from $V$ in the course of the above algorithm.
\label{def_cs}
}
\end{definition}

\noindent
As an immediate consequence of this definition we have the following
corollary.

\begin{corollary}
 If $U\subset \M$ is a finitely generated set with respect to a constructive
 involutive division $L$, then a compact involutive basis of
 ideal $Id(U)$ is minimal.
 \label{cor_cis}
\end{corollary}

\section{Minimal Involutive Bases of Polynomial Ideals}

In paper~\cite{GB96} we proposed the next algorithm {\bf InvolutiveBasis}
for computation of involutive bases of polynomial ideals. In the
algorithm the initial polynomial set $F$ is
subject, first of all, to the conventional autoreduction in line 2. Next
are two main steps which are sequentially
made:
\begin{enumerate}
\renewcommand{\theenumi}{(\roman{enumi})}
\item By the normal strategy, a non-multiplicative
prolongation $g\cdot x$ of element $g$ in the intermediate basis $G$
with the lowest $lm(g\cdot x)$ is selected in line 5. If there are several
different non-multiplicative prolongations with the same leading term, then
any of them may be selected.
\item If $h=NF_L(g\cdot x,G)\neq 0$, then $G$ is
enlarged by $h$, and the involutive autoreduction of the enlarged set
is done in line 8.
\end{enumerate}
In order to apply the criterion in line 7 for elimination
of superfluous involutive reductions and also to avoid
repeated prolongations, the auxiliary set $T$ of triples
$(g,u,P)$ is used. Here $g\in G$, and $u$
is either the lowest, with respect to the ordering $\prec$,
leading monomial in $lm(G)$ such that $g$ was produced
by non-multiplicative prolongations of $f\in G$ with $u=lm(f)$,
or $u=lm(g)$ if there is no such $f$ in $G$. Those variables in
$NM_L(g,G)$ have been chosen in line 5 are collected in set $P$.

\setcounter{cc}{00}
\vskip 0.3cm
\noindent
\h Algorithm {\bf InvolutiveBasis:}
\vskip 0.2cm
\noindent
\h {\bf Input:} $F$, a finite polynomial set
\vskip 0.0cm
\noindent
\h {\bf Output:} $G$, an involutive basis of the ideal $Id(F)$
\vskip 0.0cm
\noindent
\h \mbegin
\hln
\hh $G:=Autoreduce(F)$;\ \ $T:=\emptyset$
\hln
\hh \mfore $g\in G$\bb \mdo $T:=T\cup \{(g,lm(g),\emptyset)\}$
\hln
\hh \mwhile exist $(g,u,P)\in T$ and $x\in NM_L(lm(g),lm(G))\setminus P$
  \bb \mdo
\hln
\hhh \mchoose such $(g,u,P),x$ with the lowest
 $lm(g)\cdot x$ w.r.t. $\prec$
\hln
\hhh $T:=T\setminus \{(g,u,P)\} \cup \{(g,u,P\cup \{x\})\}$
\hln
\hhh \mif $Criterion(g\cdot x,u,T)$ is false \bb \mthen $h:=NF_L(g\cdot x,G)$
\hln
\hhhh \mif $h\neq 0$\bb \mthen $G:=Autoreduce_L(G\cup \{h\})$
\hln
\hhhhh \mif $lm(h)=lm(g\cdot x)$\bb \mthen $T:=T\cup \{(h,u,\emptyset)\}$
\hln
\hhhhh \melse $T:=T\cup \{(h,lm(h),\emptyset)\}$
\hln
\hhh $Q:=T$;\ \ $T:=\emptyset$
\hln
\hhh \mfore $g\in G$\bb \mdo
\hln
\hhhh \mif exist $(f,u,P)\in Q$ s.t. $lm(f)=lm(g)$\bb \mthen
\hln
\hhhhh \mchoose $g_1\in G$ s.t. $lm(g_1)|_Lu$
\hln
\hhhhh $T:=T\cup \{(g,lm(g_1),P)\}$
\hln
\hhhh \melse $T:=T\cup \{(g,lm(g),\emptyset)\}$
\hln
\hh \mend
\hln
\h \mend
\hln
\vskip 0.2cm
\noindent
$Criterion(g,u,T)$ is true provided that if there is $(f,v,D)\in T$
such that $lm(f)|_Llm(g)$ and $lcm(u,v) \prec lm(g)$. Correctness of this
criterion, which is just the involutive form~\cite{GB96} of
the Buchberger's chain criterion, is provided by Theorem~\ref{th_criteria}.
\vskip 0.3cm

\begin{definition}{\em Given a constructive division $L$, a finite
 involutive basis $G$ of ideal $Id(G)$ is called {\em minimal} if
 $lt(G)$ is the minimal involutive basis of the monomial ideal
 generated by $\{lt(f)\ |\ f\in Id(G)\}$.
} \label{min_ipb}
\end{definition}

\begin{theorem}
 A monic minimal involutive basis is unique.
 \label{pr_min}
\end{theorem}

\noindent
{\bf Proof}\ \ Assume for a contradiction that a polynomial ideal
 $Id(F)$ has two distinct monic minimal involutive bases $G_1$ and $G_2$.
 Their minimality means that $lm(G_1)=lm(G_2)$. Since $G_1$ and $G_2$
 are distinct there are $g_1\in G_1$ and $g_2\in G_2$ such that
 $lt(g_1)=lt(g_2)$ but $g_1\neq g_2$. Since $g_1-g_2\in Id(F)$, by
 Theorem~\ref{th_nf}, we have $NF_L(g_1-g_2,G_1)=NF_L(g_1-g_2,G_2)=0$.
 Therefore, at least one of the sets $G_1,G_2$ is not involutively
 autoreduced, and, hence, in accordance with Definition~\ref{def_inv},
 it cannot be involutive basis. \hfill{\Box}

\vskip 0.2cm
 For a globally defined involutive division, by Proposition~\ref{local_id},
 this proof, obviously, is also valid for polynomial ideals with infinite
 involutive bases. Therefore, we have the following corollary.

\begin{corollary}
 Given a globally defined involutive division, every polynomial ideal
 has the unique involutive basis.
 \label{cr_gdd}
\end{corollary}

\noindent
Thus, given a globally defined involutive division $L$, the output of
algorithm {\bf InvolutiveBasis}, in the case of its termination, is
unique for a given polynomial ideal irrespective of an ideal generating
set $F$ in the input.

However, even though the algorithm may not terminate it is still able to
compute a \Gr basis as the following proposition shows.

\begin{proposition}
 Let $L$ be a continuous involutive division and $G$ be an intermediate
 polynomial basis generated by algorithm {\bf InvolutiveBasis}. If the
 ordering $\prec$ is degree compatible, then in a finite number of steps
 $G$ becomes a \Gr basis.
\label{pr_gb}
\end{proposition}

\noindent
{\bf Proof}\ \ Let the current prolongation $g\cdot x$ is
such that $h=NF_L(g\cdot x,G)\neq 0$. Then at the second main step
of the algorithm (step (ii) as described above), the intermediate
polynomial set is enlarged by $h$. In so doing there are two
alternatives:
$$ (a)\ \ lm(h)=lm(g\cdot x);\qquad (b)\ \ lm(h)\prec lm(g\cdot x)\,.$$
In the latter case $lm(g\cdot x)$ is involutively reducible by some
$lt(f)\in lt(G)$, that is,
$lm(g)\cdot x=lm(f)\times w$. Then, by Theorem~\ref{th_nfzero}
and Corollary~\ref{c_inf}
we have the equality $NF_L(g\cdot x,G)=NF_L(S(f,g),G)$
where $S(f,g)=g\cdot x-f\times w$ is an  $S-$polynomial.

In this case, unlike the case $(a)$, the monomial ideal $Id(lm(G))$ is
changed.
Indeed, let there is a polynomial $h_1\in G$ such that $lm(h)$
is multiple of $lm(h_1)$ but not involutively multiple, that is,
$lm(h)=lm(h_1)\cdot (lm(h)/lm(h_1))$. By the
normal selection strategy, set $G$ satisfies the
condition (\ref{cond_pinv}) of partial involutivity up to the monomial
$lm(h)$ with respect to the ordering $\prec$ what implies $NF_L(h,F)=0$.

Furthermore, by Theorem~\ref{th_nf},
$NF_L(S(g_1,g_2),G)=NF(S(g_1,g_2),G)=0$
for any $S-$polynomial $S(g_1,g_2)$, ($g_1,g_2\in G$) with
$lcm(lm(g_1),lm(g_2)) \prec lm(g\cdot x)$.

It remains to prove that every $S(g_1,g_2)$ such that
$NF(S(g_1,g_2),G)\neq 0$ is computed at some step of
the algorithm. Since set $G$ is $L-$autoreduced,
monomial $u=lcm(lm(g_1),lm(g_2))$
cannot be involutively multiple of both $lm(g_1),lm(g_2)$.
Hence, by degree compatibility of the ordering
$\prec$, in a finite number of steps at least
one of $g_1,g_2$ will be non-multiplicatively prolonged to
a polynomial $g$ with $lm(g)=u$. Let $g$ be obtained by
non-multiplicative prolongations of $g_1$,
and the current prolongation be $g$ with
$u=lm(g_1)\cdot (u/lm(g_1))$.
If $u$ is involutively multiple of $lm(g_2)$ or
$lm(g_3)$ where $g_3$ is a polynomial obtained in
the course of the algorithm by non-multiplicative
prolongations of $g_2$, then we are done.

Otherwise, there is to be $\tilde{g}\in G$ such that
$u=lm(\tilde{g})=lm(g_2)\cdot (u/lm(g_2))$, and
one of the two polynomials $g,\tilde{g}$
will be constructed before the another.
Since their leading monomials coincide, the leading monomial of the
latter will be involutively reducible by the leading monomial of
the former.  \hfill{\Box}

\vskip 0.2cm
Though, by Corollary~\ref{cr_gdd}, algorithm {\bf InvolutiveBasis},
if it terminates, computes the minimal involutive basis for a globally
defined involutive division it may not be the case for arbitrary
involutive division. If we use, for instance, any of divisions in
Examples~\ref{div_T}-\ref{div_J} and \ref{div_I}, then, given a polynomial
ideal $Id(F)$, the algorithm output depends on the structure of input
generating set $F$.

\begin{example} {\em Let $F=\{x^2y-1,xy^2-1,y^4-1\}$. The lexicographical
 Janet basis for $x\succ y\succ z$ computed by algorithm
 {\bf InvolutiveBasis} is
 $$\{\ x^2y-1,x^2-1,xy^2-1,xy-1,x-1,y^4-1,y^3-1,y^2-1,y-1\ \}\,.$$
 The reduced \Gr basis $\{x-1,y-1\}$ of $Id(F)$ is also the
 minimal Janet basis.
} \label{exm_nuib}
\end{example}

\begin{proposition} If algorithm {\bf InvolutiveBasis} takes a reduced
 \Gr basis as input it produces a minimal involutive basis for a
 constructive involutive division.
 \label{gb_minib}
\end{proposition}

\noindent
{\bf Proof}\ \ Let $g\cdot x$ be a non-multiplicative
 prolongation of element $g$ in intermediate polynomial set $G$,
 and $h=NF_L(g\cdot x,G)$. We note that either $h=0$ or
 $lm(h)=lm(g\cdot x)$. Otherwise, as shown in the proof of
 Proposition~\ref{pr_gb}, $lm(h)$ would not belong to monomial
 ideal $Id(lm(G))=Id(lm(F))$. Thus, the output monomial set $lm(G)$ is
 constructed just as it would be done by applying algorithm
 {\bf MinimalInvolutiveMonomialBasis}
 to $lm(F)$. It follows that $lm(G)$ is the minimal basis of $Id(lm(F))$.
\hfill{\Box}

\vskip 0.2cm

The next algorithm constructs
a minimal involutive basis, and
generally deals with less number of intermediate polynomials than
algorithm {\bf InvolutiveBasis} causing the computational
efficiency to increase.

\begin{theorem} Let $F$ be a finite subset of $\R$ and $L$ be a
 constructive involutive division. Suppose ordering $\succ$ is degree
 compatible. Then algorithm {\bf
 MinimalInvolutiveBasis} computes a minimal involutive basis of $Id(F)$
 if this basis is finite. If $L$ is noetherian, then the basis
 is computed for any ordering.
\label{th_alg_min}
\end{theorem}

\noindent
{\bf Proof}\ \ {\em Correctness.} First of all, we recall that
correctness of the involutive criterion which is verified in
lines 14, 23 follows from Theorem~\ref{th_criteria}. As distinct from
the algorithm {\bf InvolutiveBasis} here are two disjoint subsets $T$ and
$Q$ of the triples. They are built in such a way that $lm(g)\prec lm(f)$
for any $g$ in $(g,u,P)\in T$ and $f$ in $(f,v,D)\in Q$. Let $\tilde{G}$
be a polynomial set $\{g\ |\ (g,u,P)\in Q\}$. First of all, we claim
that ideal $Id(G\cup \tilde{G})$ is an invariant of the
{\bf repeat}-loop. Indeed, it is trivially true upon
initialization. Inside the loop, if a polynomial is removed from $G$
in lines 18 and 28, then it is added to $\tilde{G}$. On the other
hand, removal of a triple from $Q$, that is, the corresponding
polynomial from $\tilde{G}$ in line 11, does not change $G$ iff
$NF_L(g,G)=0$.

Furthermore, set $T$ is handled by the lower {\bf while}-loop in lines
19-29 just as it done in algorithm {\bf InvolutiveBasis} except for
restriction in line 20 and the set contraction in lines 27-28. In the
latter case all the elements in $G$ with $lm(g) \succ lm(h)$, where
$h$ is the normal form of the current prolongation, are moved to
$\tilde{G}$. Thus, this {\bf while}-loop preserves the property
of partial involutivity up to monomial $v\prec lm(h)$ for the intermediate
set $G$, in accordance with Theorems~\ref{th_inv_cond} and \ref{th_nf}, if
there is a partially involutive set in the input of the loop. Besides,
two elements with coinciding leading terms obviously never occur in
set $\tilde{G}$.

\vskip 0.3cm
\setcounter{cc}{00}
\noindent
\h Algorithm {\bf MinimalInvolutiveBasis:}
\vskip 0.2cm
\noindent
\h {\bf Input:}  $F$, a finite polynomial set
\vskip 0.0cm
\noindent
\h {\bf Output:} $G$, the minimal involutive basis of the ideal $Id(F)$
\vskip 0.0cm
\noindent
\h \mbegin
\hln
\hh $F:=Autoreduce(F)$
\hln
\hh \mchoose $g\in F$ with the lowest $lm(g)$ w.r.t. $\prec$
\hln
\hh $T:=\{(g,lm(g),\emptyset)\}$;\ \ $Q:=\emptyset$;\ \ $G:=\{g\}$
\hln
\hh \mfore $f\in F\setminus \{g\}$\bb \mdo
\hln
\hh $Q:=Q\cup \{(f,lm(f),\emptyset )\}$
\hln
\hh \mrepeat
\hln
\hhh $h:=0$
\hln
\hhh \mwhile $Q\neq \emptyset$\bb \mand $h=0$\bb \mdo
\hln
\hhhh \mchoose $g$ in $(g,u,P)\in Q$ with the lowest $lm(g)$ w.r.t. $\prec$
\hln
\hhhh $Q:=Q\setminus \{(g,u,P)\}$
\hln
\hhhh \mif $Criterion(g,u,T)$ is false \bb \mthen $h:=NF_L(g,G)$
\hln
\hhh \mend
\hln
\hhh \mif $h\neq 0$\bb \mthen $G:=G\cup \{ h \}$
\hln
\hhhh \mif $lm(h)=lm(g)$\bb \mthen $T:=T\cup \{(h,u,P)\}$
\hln
\hhhh \melse $T:=T\cup \{(h,lm(h),\emptyset )\}$
\hln
\hhhhh \mfore  $f$ in $(f,v,D)\in T$ s.t. $lm(f) \succ lm(h)$\bb \mdo
\hln
\hhhhh $T:=T\setminus \{(f,v,D)\}$;\ \ $Q:=Q \cup \{(f,v,D)\}$;\ \
    $G:=G\setminus \{f\}$
\hln
\hhh \mwhile exist $(g,u,P)\in T$ and $x\in NM_L(g,G)\setminus P$ and,
     if $Q \neq \emptyset$,
\hln
\hhhh s.t. $lm(g\cdot x) \prec lm(f)$ for all $f$ in $(f,v,D)\in Q$\bb  \mdo
\hln
\hhhh \mchoose such $(g,u,P),x$ with the lowest $lm(g)\cdot x$ w.r.t. $\prec$
\hln
\hhhh $T:=T\setminus \{(g,u,P)\} \cup \{(g,u,P\cup \{x\})\}$
\hln
\hhhh \mif $Criterion(g\cdot x,u,T)$ is false \bb \mthen $h:=NF_L(g\cdot x,G)$
\hln
\hhhhh \mif $h\neq 0$\bb \mthen $G:=G\cup \{h\}$
\hln
\hhhhhh \mif $lm(h)=lm(g\cdot x)$ \bb \mthen $T:=T\cup \{(h,u,\emptyset )\}$
\hln
\hhhhhh \melse $T:=T\cup \{(h,lm(h),\emptyset )\}$
\hln
\hhhhhhh \mfore  $f$ in $(f,v,D)\in T$\bb with\bb $lm(f) \succ lm(h)$\bb \mdo
\hln
\hhhhhhh $T:=T\setminus \{(f,v,D)\}$;\ \ $Q:=Q \cup \{(f,v,D\})$;\ \
 $G:=G\setminus \{f\}$
\hln
\hhh \mend
\hln
\hh \muntil $Q\neq \emptyset$
\hln
\h \mend
\hln

\vskip 0.3cm
\noindent
In what follows polynomials in $\tilde{G}$, if $\tilde{G}\neq \emptyset$, are
successively selected in accordance with the normal strategy; taken
out of the set and $L-$reduced modulo $G$. The upper {\bf while}-loop in
lines 9-13 proceeds until the normal form $h$ of the selected polynomial
does not vanish. Then set $G$ is enlarged by $h$ in line 14.
The {\bf repeat}-loop terminates when set $\tilde{G}$ becomes empty in
line 11 and the lower {\bf while}-loop does not lead to
appearance of new elements in this set. It means that the output
set $G$ is an involutive basis of ideal $Id(G)=Id(F)$.

Now, by Corollary~\ref{cor_cis}, to prove minimality of the output
basis it is sufficient to show that the lower {\bf while}-loop always
ends up with $L-$autoreduced polynomial set $G$ such that $lt(G)$ is compact.
As we have already seen, this loop preserves partial involutivity.
Initially there is a single polynomial which has the minimal
leading monomial, and, therefore, its handling in the loop produces
a compact leading monomial set.

Suppose a partially involutive polynomial set $G$ with compact $lm(G)$ was
produced by the lower {\bf while}-loop, and then it is enlarged by
$h=NF_L(g,G)$ in line 14 when $G$ is partially involutive up to some
monomial $v\prec lm(g)$.

If $lm(h)=lm(g)$, then, by restriction in line 20, $lm(h)\succ lm(f)$ for
all $f\in G$. By property (d) in Definition~\ref{inv_div}, we obtain that
$NM_L(lm(f),lm(G_1=G\cup \{h\}))\subseteq NM_L(lm(f),lm(G))$ for any $f\in G$.

Let $lm(h)$ has no conventional divisors in $lm(G)$. Then, starting with
the set $G_0=Autoreduce(G_1)$, and completing $G_0$ with irreducible
non-multiplicative prolongations of its elements by the normal strategy,
we construct set $G_2\supseteq G_1$ partially involutive up to the
monomial $v$
and with compact $lm(G_2)$. If we start now with set $lm(G_1)$
and complete it, if necessary, with irreducible non-multiplicative
prolongations of its elements in order to obtain a partially involutive
set up to $v$, then we arrive at the same set $G_2$. Indeed, even in the
presence of extra intermediate elements, if $G_2\setminus G_1\neq \emptyset$,
there cannot occur reduction of an element $p\in G$ either by an
element in $G$ or by an extra element. The former reduction is
impossible by property (d) of involutive division. The latter reduction,
if it would hold, by properties (c)-(d) and by Theorem~\ref{th_nf}, would
lead to reducibility of $p$ in the earlier set $G$ when $h$
has not been added yet.

If $lm(h)$ is multiple of some element in $lm(G)$, then
continuation of processing with $G_1$ in the lower {\bf while}-loop
yields a partially involutive polynomial set up to $lm(h)$. In
doing so, $h$ is involutively reduced either to zero, or to a
polynomial which changes the monomial ideal $Id(lm(G))$, as we have
shown in the proof of Proposition~\ref{pr_gb}. Correspondingly,
$G$, after contraction in lines 27-28, is reset to the partially
involutive form with the compact leading monomial set.

In the case when $lm(h)\prec lm(g)$, the elimination which is done in
line 18 converts, apparently, the situation into one of two alternatives
we have just considered.

Thus, the {\bf repeat}-loop, if it terminates, ends up with
an involutive set $G$ with compact $lm(G)$, that is, with the minimal
involutive basis.
\vskip 0.2cm
{\em Termination}. As it shown in the proof of Proposition~\ref{pr_gb},
there may be a finite number of cases when polynomial
$g$ chosen in lines 10 or prolongation $g\cdot x$ chosen in
line 21 have reducible leading monomials. It implies finitely many
redistributions of triples between $T$ and $Q$ done in lines 18 and
28. If $Id(F)$ has the finite minimal involutive basis, and ordering $\prec$
is degree compatible, then the lower {\bf while}-loop terminates
irrespective of $Q$ is empty set or not. This follows immediately from
Propositions~\ref{pr_gb}, \ref{gb_minib} and compactness of $lm(G)$.
Since the upper {\bf while}-loop is obviously terminates, and set $Q$
is refreshed finitely many times, in a finite number of steps the
algorithm arrives at $Q=\emptyset$ in line 30.

If involutive division $L$ is noetherian then the algorithm
terminates for any ordering $\prec$ because the lower {\bf while}-loop
terminates for the same reason as the {\bf while}-loop does in
algorithm {\bf InvolutiveBasis}~\cite{GB96}.
\hfill{\Box}

\section{Conclusion}

As we noted above, algorithm {\bf MinimalInvolutiveBasis} deals, generally,
with less number of intermediate polynomials then algorithm {\bf
InvolutiveBasis}. Besides, if involutive division $L$ is not globally
defined, then we may not obtain the minimal involutive basis in the output of
the latter algorithm. But even for globally defined divisions the former
algorithm avoids the involutive autoreduction done in the latter algorithm
at every step of the intermediate set enlargement. That is why we expect
higher efficiency of algorithm {\bf MinimalInvolutiveBasis} with respect
to algorithm {\bf InvolutiveBasis} for arbitrary involutive division.

One could also construct the minimal involutive basis by
computing the reduced \Gr basis and then enlarging it by
non-multiplicative prolongations of its elements until the leading
monomial set becomes involutive. To construct the
reduced \Gr basis one can use the Buchberger algorithm or
perform the conventional autoreduction of an involutive basis
computed by algorithm {\bf InvolutiveBasis}. However,
unlike Buchberger algorithm, algorithm {\bf MinimalInvolutiveBasis}
benefits from the involutive technique, and as we have argued is
favored over the use of algorithm {\bf InvolutiveBasis} for
intermediate computation.

In paper~\cite{Schwarz} for constructing Janet bases for linear
partial differential equations one more algorithm is described which
goes back to the original computational scheme~\cite{Janet}. Its analog
in commutative algebra contains two basic subalgorithms which are
successively performed: completion of a polynomial set by non-multiplicative
prolongations of its elements until the set of leading monomials becomes
involutive or complete (see footnote at page 4); the conventional
autoreduction of the obtained set. In this case due to the second
subalgorithm the output Janet bases are minimal. However, such an
algorithmic procedure is far short of optimum from the computational
point of view. In so doing one has to perform the repeated prolongations
and deal with all the possible $S-$polynomials. In our algorithm
{\bf MinimalInvolutiveBasis} the repeated prolongations are eliminated
by storing in the triple sets $T$ and $Q$ those non-multiplicative
variables which have been used for a given polynomial. Furthermore,
the use of the involutive analogue of the Buchberger's chain criterion
allows one to cut considerably the number of computed $S-$polynomials.

The algorithms described in this paper just as Zharkov and Blinkov
algorithm can be extended to systems
of linear systems of partial differential equations~\cite{Gerdt},
and also to some classes of nonlinear systems. Being uniquely defined,
minimal involutive bases much like reduced \Gr bases can be considered
as canonical ones for polynomial and differential ideals. The
corresponding form of partial differential equation systems is just
the standard~\cite{Reid} one. By transforming a given system into
this form one can determine the dimension of the solution space
and a set of initial conditions providing the existence of a uniquely
defined and locally holomorphic solution~\cite{Janet,Riquer,Reid,Seiler}.
Involutive algorithmic ideas may be also rather fruitful in constructing
the canonical bases for finitely generated ideals in free Lie algebras
and superalgebras~\cite{GK96}.

\section{Acknowledgements}
 The authors are grateful to Joris van der Hoeven for useful comments and
 suggestions on Division II, and to Joachim Apel for important remarks.
 This work was supported in part by the RFBR grant No. 96-15-96030.

\end{document}